\documentclass[12pt]{amsart}
\usepackage{amscd}
\usepackage{amssymb}
\usepackage{amsfonts}
\usepackage[matrix,arrow,curve]{xy}

\newcommand{\CH}{\operatorname{CH}}
\newcommand{\tdiv}{\operatorname{div}}

\newcommand{\Alb}{\operatorname{Alb}}
\newcommand{\Pic}{\operatorname{Pic}}
\newcommand{\lin}{\operatorname{lin}}
\newcommand{\red}{\operatorname{red}}
\newcommand{\sm}{\operatorname{sm}}
\newcommand{\NS}{\operatorname{NS}}

\newcommand{\tIm}{\operatorname{Im}}
\newtheorem{thm}{Theorem}[section]
\newtheorem{lemma}[thm]{Lemma}

\newtheorem{prop}[thm]{Proposition}

\theoremstyle{definition}
\newtheorem{df}[thm]{Definition}
\newtheorem{theosd}[thm]{Theorem}
\newtheorem{Propsd}[thm]{Proposition}

\newtheorem{rem}[thm]{Remark}

\numberwithin{equation}{section}

\begin{document}

\title{Families of disjoint divisors on varieties}

\author{Fedor Bogomolov, Alena Pirutka and Aaron Silberstein}

\address{Fedor Bogomolov,  Courant Institute of Mathematical Sciences, New York University, 251 Mercer Street, New York, NY 10012-1185 and  National Research University Higher School of Economics, Russian Federation AG Laboratory, HSE, 7 Vavilova Str., Moscow, Russia, 117312}
\email{bogomolo@cims.nyu.edu}

\address{Alena Pirutka, Courant Institute of Mathematical Sciences, New York University, 251 Mercer Street, New York, NY 10012-1185 }

\email{pirutka@cims.nyu.edu}

\address{Aaron Silberstein,  The University of Chicago, 5734 S.~University Ave., Room 208C, Chicago, IL 60637}

\email{asilbers@math.uchicago.edu}

\begin{abstract}
Following the work of Totaro and Pereira, we study sufficient conditions under which collections of pairwise-disjoint divisors on a variety over an algebraically closed field are contained in the fibers of a morphism to a curve. We prove that $\rho_{w}(X) + 1$ pairwise-disjoint, connected divisors suffices for proper, normal varieties $X$, where $\rho_{w}(X)$ is a modification of the N\'eron-Severi rank of $X$ (they agree when $X$ is projective and smooth).  We then prove a strong counterexample in the affine case: if $X$ is quasi-affine and of dimension $\geq 2$ over a \textit{countable}, algebraically-closed field $k$, then there exists a (countable) collection of pairwise-disjoint divisors which cover the $k$-points of $X$, so that for any non-constant morphism from $X$ to a curve, at most finitely many are contained in the fibers thereof. We show, however, that an uncountable collection of pairwise-disjoint, connected divisors in \textit{any} normal variety over an algebraically-closed field must be contained in the fibers of a morphism to a curve.
\vspace{0.5cm}

\end{abstract}

\maketitle

\section{Introduction}\label{intro}

The goal of this note is to give a set-theoretic condition under which collections of pairwise-disjoint divisors on varieties over an algebraically-closed field are contained in the fibers of a single morphism to a curve.  We first adapt the methods of B.~Totaro \cite{TotaroDisjDiv} and J.~Pereira \cite{PereiraFol} to produce a stronger bound in the projective, smooth case in characteristic zero, and we generalize these  results to normal, proper varieties in all characteristics.
 
We obtain the following result, which generalizes  the theorems of Totaro and Pereira, \textit{loc.~cit.}, to normal, proper varieties of arbitrary characteristic.  Here, $\rho_{w}(X)$ is an invariant of the variety $X$, equal to the N\'eron-Severi rank when $X$ is smooth and projective, and finite in all cases :

\theosd\label{DDP}{\it Let $X$ be a normal, proper, integral variety defined over an algebraically closed field $k$.  Let $\{D_i\}_{i\in I}$ be a collection of pairwise-disjoint, reduced, codimension-one, connected subvarieties of $X$.  Assume that $\# I\geq \rho_w(X)+1$. Then there is a smooth, projective curve $C$ and a  surjective morphism $f:X\to C$ with connected fibers such that for any $i \in I,$ the divisor $D_i$ is contained in a fiber of $f$.  Furthermore, there is a set $\Sigma \subseteq I$ so that $\#(I \setminus\Sigma) \leq \rho_{w}(X) - 2$ and for each $i \in \Sigma, D_{i}$ is equal (set-theoretically) to a fiber of $f$.  \\ }

In Pereira and Totaro's approaches --- which work only in the smooth case --- $\#I$ must be at least $\rho_{w}(X) + 2$; our extra saving comes from an extra application of the Hodge index theorem.
  
In the affine case we have the following explicit counterexample.\\

\theosd\label{DDA}{\it Let $\mathbb{A}^{2}_{k}$ be the affine plane over a countable, algebraically-closed field $k$. Then there is a countable family   $\{D_i\}_{i\in I}$  of  integral, Zariski-closed, codimension-$1$ subvarieties of $X$, such that:
\begin{itemize} \item[$\bullet$] The divisors $D_i$  are pairwise disjoint and their $k$-points cover $\mathbb{A}^{2}$, i.e. $\mathbb{A}^{2}(k)=\bigcup_{i\in I} D_i(k)$;
\item[$\bullet$]  For any non-constant morphism $f: \mathbb{A}^{2}\rightarrow C$ to a curve, at most finitely many of the $D_{i}$ are contained in fibers of $f$.\\
 \end{itemize}}
As a corollary, we easily deduce:
\theosd\label{DDArb}{\it Let $X$ be a quasi-affine variety over a countable, algebraically-closed field $k$. Then there is a countable family $\{D_i\}_{i\in I}$  of  connected, Zariski-closed, codimension-$1$ subvarieties of $X$, such that:
\begin{itemize} \item[$\bullet$] The divisors $D_i$  are pairwise disjoint and their $k$-points cover $X$, i.e., $X(k)=\bigcup_{i\in I} D_i(k)$;
\item[$\bullet$]  For any non-constant morphism $f: X\rightarrow C$, at most finitely many of the $D_{i}$ are contained in fibers of $f$.\\
 \end{itemize}}
We can salvage this counterexample if $I$ is \textit{uncountable}:
 \theosd\label{DDAU}{\it Let $X$ be \textit{any} normal variety over an algebraically-closed field $k$, and let $\{D_{i}\}_{i\in I}$ be an \textit{uncountable} collection of pairwise-disjoint, reduced, codimension-one, connected subvarieties of $X$.  Then there is a normal curve $C$ and a non-constant morphism $f: X\rightarrow C$ with connected fibers so that each $D_{i}$ is contained in a fiber of $f$. \\ \\
 }
In particular, if a set of divisors covers the $k$-points of the variety when $k$ is uncountable, there are uncountably many divisors.% each divisor is defined over a finitely-generated subfield of $k$, and $k$ must have uncountable transcendence degree over $\overline{\mathbb{Q}}$, so the theorem holds for partitions of the $k$-points of $X$ by divisors.

This paper provides a tool to approach the third author (A.S.)'s program of Geometric Reconstruction \cite{AST} in the first author (F.B.)'s Program of Birational Anabelian Geometry.  In Bogomolov's program, we take a field $K$ which is the function field of an algebraic variety $X$ of dimension $\geq 2$ defined over an algebraically-closed field, and the goal is to reconstruct $K$ from its absolute Galois group $G_{K}$.  In Geometric Reconstruction, the goal is to reconstruct individual \textit{varieties}   with a given function field $K$ as group-theoretically defined objects in $G_{K}$.  The results in this paper are crucial  in an upcoming paper by the third author in proving geometric reconstruction for function fields $K$ of transcendence degree $2$ over $\overline{\mathbb{Q}}$ from the maximal, $2$-step nilpotent, pro-$\ell$ quotient of $G_{K}$.

\section{Disjoint divisors on proper varieties}

In this section we prove Theorem \ref{DDP}, following the proof of Totaro \cite{TotaroDisjDiv}, but with two additional arguments.  Totaro and Pereira prove the theorem in characteristic zero, for smooth, projective varieties.  First, in order to generalize to characteristic $p$, we reduce the theorem for $X$ normal and projective to the case of a general surface which is an intersection of hyperplane sections; this allows us to use resolution of singularities of surfaces, in arbitrary characteristic.  By appealing to the Hodge index theorem we reduce the number of pairwise-disjoint divisors to the theoretical minimum.

\subsection{Divisors and Albanese varieties}

We begin with the notion of divisor class group we will use throughout the paper. Recall that for a variety $Y$ over a field $k$ we denote by $Z^1(Y)$  the group of Weil prime divisors on $X$ --- that is, finite, linear combinations of closed, integral, codimension-$1$ subvarieties and   by $CH^1(Y)$ the quotient of $Z^1(Y)$ by linear equivalence. 

 Let $X$ be a  normal projective  integral variety defined over an algebraically closed field $k$ and $U\subset X$ be the smooth locus of $X$.  Fix a prime $\ell$ not equal to the characteristic of $k$.  Since $X$ is normal, the singular locus is codimension $\geq 2$ and the restriction map $\CH^1(X)\to \CH^1(U)$ is an isomorphism.
  Since $U$ is smooth, we have the cycle class map $CH^1(U)\to  H^2_{\acute{e}t}(U, \mathbb Q_{\ell}(1))$ \cite[2.1.1]{SGA45}.   We  then denote $CH^1(X)_{hom}\subset CH^1(X)$ the group of cycles homologically equivalent to zero as the kernel of the composition  $$CH^1(X)\stackrel{\sim}{\to} CH^1(U)\to  H^2_{\acute{e}t}(U, \mathbb Q_{\ell}(1))$$ for $\ell\neq char\,k$.  
  If $X$ is just normal and proper, we can use Chow's lemma \cite[Th.~5.6.1]{EGA2} to find a \textbf{projective, normal modification} $\pi: \tilde{X}\rightarrow X$ --- where $\pi$ is a projective, surjective, birational morphism, and $\tilde{X}$ is a normal, projective variety (we can assume normality because normalization is a projective morphism). For a normal, projective variety $Y,$ define the group $B^{1}(Y):= (\CH^{1}(Y)/\CH^{1}(Y)_{\hom})\otimes \mathbb{Q}$

 Note that if $X$ is smooth, then  algebraic equivalence (see \cite[10.3]{IntersectionTheory})  and numerical equivalence tensored with $\mathbb Q$ coincide  for codimension $1$ cycles (\cite[6.3]{Kleiman}), also numerical and homological equivalences  tensored with $\mathbb Q$  coincide for codimension $1$ cycles \cite[3.4.6.1]{Andre}. In particular, $B^{1}(Y) = \NS(Y)\otimes \mathbb{Q}$.  % so that $\rho_w(X)=\rho(X)$ is the rank of the N\'eron-Severi group of $X$. In particular, $\rho_w(X)$ is independent of $\ell$.

 \lemma{ Let $X$ be a  normal projective  integral variety defined over an algebraically closed field $k$ and $U\subset X$ be the smooth locus of $X$. Then:
 \begin{itemize}\item[(i)] the  algebraic and homological equivalence tensored over $\mathbb Q$ coincide for codimension $1$ cycles on $X$; 
 \item[(ii)] $\dim_\mathbb{Q}B^1(X)< \infty$, and this number is independent of the $\ell$ chosen in the definition.
 \end{itemize}}
 \proof{Note that if a codimension $1$ cycle $\alpha\in Z^1(U)$ is algebraically equivalent to $0$, then it is homologically equivalent to $0$. Hence, for (i), it is enough to show that if the class of $\alpha$ is $0$, then $\alpha$ is algebraically equivalent to zero. Let $f:\tilde X\to X$ be a smooth projective alteration of degree $d$, such that for $Z=X\setminus U$,  $\tilde Z=f^{-1}(Z)$ is a simple normal crossings divisor \cite{Il}. Let $\tilde U$ be the inverse image of $U$ and $U_0\subset U$ an open subset, such that the induced morphism $\tilde U_0=f^{-1}(U_0)\to U_0$ is finite of degree $d$. Note that one can assume that the complement of $U_0$ in $U$ is of codimension at least $2$: indeed,  this follows from the fact that in  the Stein factorisation $\tilde X\to Y\to  X$  of $f$, with $\tilde X\to Y$ birational and $Y\to X$ finite of degree $d$, the map $\tilde X\to Y$ is an isomorphism over any point of codimension $1$, since $X$, and so $Y$, is normal.
 
 Assume now that $\alpha$ is homologically equivalent to $0$ on $X$. Since $\tilde Z$ is the simple normal crossings divisor, we deduce that for a codimension $1$ cycle $\beta$ supported on $\tilde Z$, one has that $f^*\alpha+\beta$ is homologically equivalent to $0$ on $\tilde X$. Since $X$ is smooth, the discussion above the lemma shows that $N(f^*\alpha+\beta)$ is algebraically equivalent to $0$ on $\tilde X$, for some integer $N$. Hence, its restriction $Nf^*\alpha$ to $\tilde U_0$ is algebraically equivalent to $0$. Since $\tilde U_0\to U_0$ is finite of degree $d$, we deduce that $N\alpha$ is algebraically equivalent to $0$ on $U_0$, hence on $U$, as the complement of $U_0$ in $U$ is of codimension at least two, so that we obtain (i).
 
 For (ii), the independence of $\ell$ follows from (i). For the finiteness, it is enough to show that there is no infinite collection of divisors in $X$ with $\mathbb Q$-linearly independent classes in $H^2_{\acute{e}t}(U_0, \mathbb Q_{\ell}(1))$. Via the trace map  $H^2_{\acute{e}t}(\tilde U_0, \mathbb Q_{\ell}(1))\to  H^2_{\acute{e}t}(U_0, \mathbb Q_{\ell}(1))$  (\cite[Expos\'e IX (5.1.4)]{SGA4}) it is enough to establish the same property for $\tilde U_0$, which follows from the fact that the N\'eron-Severi group of the smooth variety $\tilde X$ is finitely generated.   \qed\\}

 The lemma above allows us to make the following definition, independently of $\ell$:

\df{The {\it Weil divisor rank} $\rho_w(X)$ of $X$ is the minimum dimension of the $\mathbb{Q}$-vector space $B^{1}(\tilde{X})$ over all projective, normal modifications $\pi: \tilde{X}\rightarrow X$ of $X$.\\ }

 In what follows we will fix $\pi: \tilde{X}\rightarrow X$ a projective, normal modification for which the dimension of the group $B^{1}(\tilde{X})$ is minimized.\\

%We will use two pullback morphisms consistently.  The first is, 

Let $X$ be normal and projective. Given a linear section $\iota: Y\hookrightarrow X $ the intersection map $\iota^{*}: \CH^{1}(X)\rightarrow \CH^{1}(Y)$, proven to be well-defined in \cite[Prop.~2.6]{IntersectionTheory}.  
%The second is more subtle:

 \Propsd\label{funct}{\it  Let $X$ be normal and projective. Let
 $\iota: T'\rightarrow X$ be an intersection of hyperplane sections of some projective embedding of $X$, smooth on the intersection with the smooth locus of $X$ (such sections are generic by \cite[Theorem 1]{Se}), of dimension $2$. 
Let $\beta: T\rightarrow T'$ be a resolution of singularities which is an isomorphism on the smooth locus of $T'$ and for which the exceptional divisors are simple normal crossing and let 
  $\phi: T\rightarrow X$ be  the composition $\iota\circ \beta$. Then the composition $p_T:\beta^{-1}\circ \iota^{*}$: induces a pullback homomorphism $p_{T}: B^{1}(X)\rightarrow B^{1}(T)$.}
 \proof{
Let $\Gamma$ be the subgroup of $\CH^{1}(T)$ generated by the exceptional divisors of $\beta$, and let $\eta: V\rightarrow T$ be the inclusion of the inverse image in $T$ of the smooth locus of $T'$.  Each exceptional divisor is in the kernel of the flat pullback morphism $\eta^{*}: \CH^{1}(T)\rightarrow \CH^{1}(V)$.    The image of $\Gamma\otimes \mathbb{Q}_{\ell}$ under the cycle map in $H^{2}_{\acute{e}t}(T, \mathbb{Q}_{\ell}(1))$ is exactly the kernel of the restriction morphism to $H^{2}_{\acute{e}t}(V, \mathbb{Q}_{\ell}(1))$ by inductive use of the Gysin sequence \cite[Cor.~16.2]{MilneEC}, so by functoriality of the cycle map, if  $\alpha\in p_{T}(\CH^{1}_{\hom}(X))$, then $\alpha$ is a linear combination of exceptional divisors on $T$. If $H$ is an ample divisor on $T'$ not passing  by $T\setminus \beta(V)$ we then deduce that  $\beta^{-1}H\cdot\alpha=0$, hence $\alpha$ is numerically equivalent to $0$ on $T$ using the Hodge index theorem for surfaces, and $p_{T}$ defines a homomorphism from $B^{1}(X)\rightarrow B^{1}(T)$. \qed\\
}

We now recall some facts about Albanese varieties.  
 
Let $Y$ be a variety defined over an algebraically closed field $k$, $Y_{i}$ its irreducible components, and let $x_{0, i}$ be a smooth point of each $Y_{i}$. We say a rational map (resp.~morphism) $f:Y\to A$ with $A$ an abelian variety is {\it admissible} if $f$ is defined at each $x_{0, i}$ and $f(x_{0, i})=0$. Following R.~Ghorpade and G.~Lachaud \cite[Section 9]{GL}, we call an Albanese-Weil variety $\Alb_w(Y)$ (resp., an Albanese-Serre variety $\Alb_s(Y)$) of $Y$ an abelian variety $A$ over $k$ with an admissible rational map $f$  (resp., morphism) from $Y$ to $A$, such that the following universal property holds: for any admissible rational map $g$  (resp., morphism) from $Y$ to an abelian variety $B$ there is a homomorphism of abelian varieties $\tilde{g} :A\to B$ such that $g=\tilde{g}\circ f$. We have the following properties:
 \begin{enumerate}
 \item The variety $\Alb_{w} (Y)$ and the universal rational map $Y\rightarrow \Alb_{w} (Y)$ exist, are independent of the choice of $x_{0, i}$ up to a translation, and $\Alb_{w}(Y) = \prod_{i} \Alb_{w}(Y_{i})$.
 \item If $Y$ is normal, the variety $\Alb_s(Y)$ exists, and is dual to the reduced Picard variety $(\Pic^0_{Y/k})_{\red}$ \cite[Paragraph after Example 9.2]{GL}.
 \item If $Y$ is smooth, the variety $\Alb_s(Y)$ coincides with $\Alb_wY$ and for $Y$ normal, there is a canonical  surjective map $\nu: \Alb_w(Y)\to \Alb_s(Y)$ with connected kernel \cite[Prop.~9.1]{GL}.
 \item A birational morphism $Y\to X$ of varieties induces an isomorphism $\Alb_{w}(Y) \to \Alb_{w}(X)$ (this follows straight from the definition), so a resolution of singularities $Y\to X$ induces an isomorphism $\Alb_{s}(Y)\to\Alb_{w}(X)$.
 \end{enumerate}
 
 We  need the following Lefschetz-type property \cite[Prop.~9.4]{GL}:
 \Propsd\label{Ltz}{\it Let $X\hookrightarrow \mathbb{P}^N$ be an embedding of $X$ into a projective space. If $i:Y\hookrightarrow X$ is a general linear section of $X$ of dimension $d\geq 2$, the canonical map $i_*:\Alb_w(Y)\to \Alb_w(X)$ induced by $i$ is a purely inseparable isogeny.}  

\vspace{0.7cm}

\subsection{The torsion case}\label{torsioncase}
$\,$\\

We start with the following easy lemma: 

\lemma\label{lineq}{
Let $\Delta_1$ and $\Delta_2$ be two effective, disjoint divisors on a proper, normal 
variety $Y$ over a field $k$, and suppose that  
\begin{equation}
\Delta_1 -\Delta_2\sim_{\lin} 0.
\end{equation}
  Then there exists  a morphism $f: Y\rightarrow \mathbb P^1_k$ such that
$\Delta_1=f^{-1}(0)$ and $\Delta_2=f^{-1}(\infty)$.
}
\proof{
By assumption, there exists a rational function $g$ on $Y$ such that $\Delta_1 -\Delta_2=\tdiv(g)$. Then we define a map  $f: Y\rightarrow \mathbb P^1_k$ 
by $f(x)=[g(x):1]$ if $x$ is not in the support of $\Delta_2$ and $[1: 0]$ (equivalently, $[1: g(x)]$ if $x$  is not in the support of $\Delta_1$) otherwise.  Since the divisors $\Delta_1$ and $\Delta_2$ are disjoint, we get a well-defined map as required. \qed\\
}

Let now $X$ be as in Theorem \ref{DDP}. We assume first that $X$ is projective.
Since   $\# I\geq \rho_w(X)+1$, there  is  a subset $J\subseteq I$ and a nontrivial linear combination  $D=\sum_{j\in J} \lambda_jD_j\in Z^1(X),$ where $\lambda_j\in \mathbb{Z}$ and $D \in \CH^{1}_{\hom}(X)$.

\prop\label{tor}{If there exists $N>0$ such that $ND\sim_{\lin} 0$, then there is a  surjective morphism $f:X \to \mathbb P^1_k$ such that for any $j\in J$ the divisor $D_j$ is contained in a fiber of $f$.}
\proof{It suffices to write $ND=\Delta_1 -\Delta_2$ as a difference of two effective (and disjoint) divisors and apply Lemma \ref{lineq}.   We obtain a map   $f: X\rightarrow \mathbb P^1_k$  satisfying the required properties: note that for each $i\in I\setminus J$ the divisor $D_i$ is a subset of a fiber of $f$:  otherwise,  for each $j \in J$, there would exist some $i \in I$ so that $D_j$ would intersect $D_i$, contradicting disjointness. \qed\\
 }
 
% If $D$ is not a torsion point in $CH^1(X)_{hom}$, then the same holds for a general linear section of $X$:
 To handle the non-torsion case, we need that non-torsion elements of $\CH^{1}$ specialize under generic hyperplane sections to non-torsion elements.
 
 \prop\label{nontorsection}{Let $X\hookrightarrow \mathbb{P}^N$ be an embedding of $X$ into a projective space. If $D$ is a non-torsion element of $\CH^1(X)$ then for a general linear section $\tau:Y\subset X$ of dimension $d\geq 2$ the restriction $\tau^*D$ given by intersection of $D$ with $Y$ to $\CH^1(Y)$ is also a non-torsion element.}
\proof{By induction, it suffices to prove the theorem for general $Y$ of codimension $1$.  We may assume $Y$ is normal \cite[Theorem 7]{Se}, and that $Y$ contains no irreducible component of $D$, as $Y$ is basepoint-free.  Suppose there is an integer $N_Y$ and a function $f_Y\in k(Y)$ such that $N_Y\tau^*D=\tdiv(f_Y)$ in $Z^1(Y)$. We can lift the function $f_Y$ to an element  $F_Y\in\mathcal O_{X,Y}\subset k(X)$. Define $D' := N_YD-\tdiv(F_Y)=\sum a_i Z_i$ with $a_i\in \mathbb Z$ and $Z_i$  irreducible components of  $D'$, that are included in $X\setminus Y$ by construction. Since $Y$ is ample, $Y$ intersects every proper codimension-$1$ subvariety, so $X$ has  no proper codimension-$1$ subvariety contained in $X\setminus Y$, so that we have $\tdiv(F_Y)-N_YD=0$ in $\CH^1(X)$, contradicting our assumption on  $D$.  \qed\\}

\subsection{The Hodge Index Theorem and the General Case}\label{generalcase}

Let $\tilde{S}$ be a projective resolution of singularities of a generic, normal, linear surface section $S\subseteq X$, whose smooth locus is exactly the intersection of the smooth locus of $X$ with $S$.
As a simple linear-algebraic corollary of the Hodge index theorem, we have:
\Propsd\label{ortho}{{\it Let $H = \{H_{j}\}_{j\in J}\subseteq B^{1}(\tilde{S})$ be an orthogonal subset of nonzero elements -- that is, $H_{j}\cdot H_{j'} = 0$ for each $j \neq j'$.  Assume furthermore that $H$ is contained in a subspace $V\subseteq B^{1}(\tilde{S})$ of dimension $d$ for which there exists $M \in V$ so that $M\cdot M > 0$. Define: 
\[
J_{+} (\text{resp.,}\,J_{-}, J_0) := \left\{j \in J \mid H_{j}\cdot H_{j} > \text{(resp.,}\,<, =)\,0\right\}.
\]
Then:
\begin{enumerate}
\item $J_{+}\cup J_{-}$ is a linearly independent set, and $\#J_{+}\leq 1$ and $\#J_{-} \leq d - 1$.
\item If $(J_{+} \cup J_{0}) \geq 2$ then $\#J_{-} \leq d - 2, \#J_{+} = 0$ and $\#J_{0} \geq \#J - (d - 2)$.
\item The span of $J_{0}$ is at most one-dimensional.
\end{enumerate}}}

The pullback $\{\tilde{D}_{i}\}_{i \in I}$ of $\{D_{i}\}_{i \in I}$ to $B^{1}(\tilde{S})$ is an orthogonal set of nonzero elements, and is contained in the image of $B^{1}(X)$, by Proposition \ref{funct}. The restriction of a general ample divisor to $\tilde{S}$ is likewise ample, so the image of $B^{1}(X)$ in $B^{1}(\tilde{S})$ is a subspace of dimension $\leq \rho_{w}(X)$ which contains an element of positive self-intersection.  By Proposition \ref{ortho}, $\#(J_{+} \cup J_{-}) \geq 2$, so in fact $\#J_{0} \geq 3$.  Let $i, j, t$ be distinct elements of $J_{0}$, and let $F := a\tilde{D}_{i} - b\tilde{D}_{j}$ be an integral linear combination in $\CH^{1}(\tilde{S})_{\hom}$.  Let $I' := I \setminus \{i, j\}$.
\prop\label{nsur}{If $F$ is not torsion, then for each $l \in I'$, the map $$\Alb_w D_{l}\rightarrow \Alb_w X$$ is not surjective.

\proof{ Fix a projective embedding $X\subset \mathbb{P}^N$.  By Proposition \ref{Ltz}, if $\tau:S\subset X$ is a general linear section of $X$ of dimension $2$, then $S$ is normal, and we have an isogeny $\Alb_wS\to \Alb_w X $. By Proposition  \ref{nontorsection}, the restriction of $F'$ of $F$ to  $\CH^1(S)_{\hom}$ is not a torsion element.  Let $\nu: \tilde S\to S$ be a resolution of singularities. Since $S$ is normal, we may assume that $\nu$ is an isomorphism over the smooth locus $S^{\sm},$ which contains all codimension one points of $S$. Let $\tilde{D}_{l}$ be a union of normalizations of the components of $\nu^{-1}(D_{l}\cap S)$ (that is, the inverse image of the intersection of $D_{l}$ intersected with $S$).  Since the support of $F$ is disjoint from $D_l$, the line bundle defined by $F$ becomes trivial on $\tilde{D}_{l}$ and we conclude that $\tilde{D}$ is a non-torsion element in the kernel of the map $\Pic^0(\tilde S)\to \Pic^0(\tilde{D}_{l})$. By duality, we obtain that the map $\Alb_s(\tilde{D}_{l})\to \Alb_s(\tilde S)$ is not surjective.  The maps $\Alb_s (\tilde S)\to \Alb_w (S)\to \Alb_w(X)$  are isogenies, so $\Alb_w(\tilde{D}_{l})\to \Alb_w(\tilde{X})$ cannot be surjective. \qed \\ \\}

{\it Proof of theorem \ref{DDP}.\\}

In Proposition \ref{tor}, we established the result if $F$ is a torsion element in $\CH^{1}(X)$.

Consider now the general case.  Let $F \in \CH^{1}(X)$ be non-torsion. Let $i, j, t$ be distinct elements of $J_{0}$ as before.  By Proposition \ref{nsur},  the map $\Alb_w(\tilde{D_{t}})\to \Alb_w(X)$ is not surjective.

By the universal property of the Albanese variety we have the following commutative diagram 
\begin{center}\mbox{
\xymatrix{ \tilde{D}_{t}\ar@{-->}[r]\ar[d] & \Alb_w(\tilde{D}_{t})\ar[d] & \\
X\ar@{-->}[r] & \Alb_w(X)\ar[r] & \Alb_w(X)/\Alb_w(\tilde{D}_{t}).
}}\end{center}
showing that $\tilde{D}_{t}$ is contracted by the composite rational map $g:X\dashrightarrow \Alb_w(X)/\Alb_w(\tilde{D}_{t}).$  

Let us show that the image of $g$ is a curve:
\begin{enumerate}
\item  We have $\dim\,\tIm(g)>0$ since the image of $X$ in $\Alb_w(X)$ generates the abelian variety $\Alb_w(X)$ by the universal property and  $\Alb_w(D_{t})\to \Alb_w(X)$ is not surjective.
\item If $\dim\,\tIm(g)$ were greater than $1$, the image of the dimension of $g$ restricted to $\tilde{S}$ would also have image of dimension $2$.  The morphism $g$ restricted to $\tilde{S}$ is a \textit{regular} map, because it factors through $\Alb_{s}(\tilde{S})/\Alb_{w}(\tilde{D}_{t})$, and the map from $\tilde{S}\rightarrow \Alb_{s}(\tilde{S})$ is defined everywhere.  Any effective divisor on $\tilde{S}$ contracted by $g$ would need to have negative self-intersection \cite[Remark after Theorem 16.2]{Artin}.  But $D_{t}$ is contracted and has self-intersection $0$, so the image has to be of dimension $<2$. 
\end{enumerate}

We then see immediately that $\tilde{D}_{j}$ is also  contracted and should be (a multiple of) a fiber of $g$ restricted to $\tilde{S}$, since its self-intersection is zero.  Therefore, $D_{j}$ and $D_{t}$ are (multiples of) fibers of the rational map to $\Alb_{w}(X)/\Alb_{w}(D_{t})$; as they are disjoint, and $X$ is normal, $g$ is in fact a \textit{regular} map --- that is, it is defined everywhere.

Therefore, the image of $ g:X\to \Alb_w(X)/\Alb_w(D_{t})$ is a curve $C'$, and all the divisors $D_{l}$ for $l \in I'$ are contained in fibers. If now $X\stackrel{f}{\to} C\to C'$ is the Stein factorization of $g, C$ is a normal curve, in which the $D_{l}$ are contained in fibers for $l \in I'$.  By \cite[Fact 1.1d]{ArtinWinters}, $D_{t}$ is in fact a multiple of an entire fiber of $f$; as $D_{i}$ and $D_{j}$ are disjoint from $D_{t}$, both $D_{i}$ and $D_{j}$ are contained in fibers, and the $D_{\sigma}$ which are (multiples of) fibers of $f$ are exactly those for which $\sigma \in J_{0}$.  We may thus set $J$ to be $J_{0}$.

Finally, if $X$ is normal and proper, consider $\{\pi^{-1}(D_{i})\}_{i \in I}$, where  $\pi:\tilde X\to X$ is a birational morphism given by Chow lemma, with $\tilde X$ projective. By Zariski's main theorem, each $\pi^{-1}(D_{i})$ is connected.  Then there exists a function $\tilde{f}: \tilde{X}\rightarrow C$ for which each of the $\pi^{-1}(D_{i})$ is contained in a fiber, and for at least two (in fact, three) $m_{1}, m_{2}, m_{3} \in I, \pi^{-1}(D_{m_{i}})$ is a (multiple of) a fiber. As $X$ is normal, by Zariski's Main Theorem, to check that $\tilde{f}$ factors through a function to $X$ we must merely check that none of the divisors that $\pi$ contracts intersects $\pi^{-1}(D_{m_{1}})$ and $\pi^{-1}(D_{m_{2}})$.  But if a divisor that $\pi$ contracts intersected both of them, the $D_{m_{i}}$ would not be disjoint. \qed \\

\section{Disjoint divisors on affine varieties}

In this section we prove Theorem \ref{DDA} and Theorem \ref{DDArb}.  Let $y, z$ be the coordinates of $\mathbb A^2_k$. Consider the following family, constructed recursively:

\begin{enumerate}
\item Define $d_{0} := 1, f_0(y, z):=z^{d_{0}}$, and $D_0 := V(f_{0}),$ the zero locus of $f_{0}$.
\item Define $d_{1} := 2, f_1(y,z):=yz^d_{1}+1$; and $D_1:= V(f_{1}).$
\item  Let $P_{2} \in \mathbb A^2_k\setminus (D_0 \cup D_1)$ and define $$a_2:=-f_1(P)/f_{0}(P)^{2d_1-1}.$$ Since $P_{2}\notin D_1,$ $a_{2} \neq 0$.
  Define  $d_{2}:= 2d_{1}$; $$f_2(y,z):=a_2f_{0}(y,z)^{d_{2}-1}+f_{1}(y,z)^d_{1};$$ and $D_2:= V(f_{2})$.  Note that in each case, $d_{i} = \deg_{z} f_{i},$ the degree of $f_{i}$ as a polynomial in $z$.
  \item Let now $n > 2$; we give a recursive definition of $f_n,$ given $f_{i}$ when $i<n$.  We define $D_{i} := V(f_{i})$ and $d_{i} := \deg_{z}(f_{i})$. Define now $d_{n} := \sum \limits_{i=1}^{n-1} d_i$.  Let $P_{n} \in \mathbb{A}^2_k\setminus \bigcup\limits_{i=0}^{n-1}D_i;$ define $$a_n:=-\prod\limits_{i=1}^{n-1} f_i(P)/f_{0}(P)^{d_n-1}.$$ Again, $a_n\neq 0$. Define $$f_n(y,z):=a_nf_{0}(y,z)^{d_n-1}+\prod \limits_{i=1}^{n-1} f_i(y,z)$$ and $D_n:= V(f_{n})$; by construction, $d_{n} = \deg_{z} f$.\\
\end{enumerate}
If we enumerate the $k$-points of $\mathbb{A}^{2}_{k}$, we may choose our $P_{i}$ so that the $\bigcup_{i}D_{i}(k) = \mathbb{A}^{2}(k)$.  By construction, the $D_{i}$ are pairwise-disjoint: the radical of the ideal generated by $f_{n}$ and $f_{i}$ for $0 < i < n$ contains $f_{0}$ and $f_{i}$; and by construction, $f_{i}$ and $f_{0}$ have no common zeroes.

We now prove that each $f_{n}$ is irreducible, and that no infinite subset of the $D_{i}$ could be contained in the fibers of a non-constant morphism.

To prove that $f_{n}$ is irreducible, we will change coordinates.  We view $\mathbb{P}^{1}_{k} = \mathbb{A}^{1}_{k}\cup \{\infty\}$.  Let $\overline{X} := \mathbb{P}^{1}_{k}\times \mathbb{P}^{1}_{k}$, and use the coordinates $y, z$ to embed $\mathbb{A}^{2}_{k}$ as an open subset $$\mathbb{A}^{2}_{k}\stackrel{\sim}{\to} \mathbb{A}^{1}_{k}\times \mathbb{A}^{1}_{k} \hookrightarrow \overline{X}$$; call this open subset $X_{1}$.   Let $\overline{D_{i}}$ be the Zariski closure of $D_{i}$ in $\overline{X}$.  Define $X_{2}\subseteq \overline{X}$ as $\mathbb{A}^{1}_{k}\times (\mathbb{P}^{1}_{k} \setminus \{0\});$ this is isomorphic to $\mathbb{A}^{2}$ with coordinates $y, x:= \frac{1}{z}$. Let $D_{i}':= \overline{D_{i}} \cap X_{2}$.

The defining ideals for $D_{i}'$ are generated by:
\begin{enumerate}
\item $f_1'(y,x)=y+x^2$;
\item $f_2'(y,x)=a_2x+(y+x^2)^2$;
\item $f_n'(y,x)=a_nx+\prod\limits_{i=1}^{n-1} f_i(y,x)$.
\end{enumerate}
By Eisenstein's criterion, applied to the ring $k(x)[y]$, all the polynomials $f_n'$ are irreducible.  Therefore, to check whether the $f_{i}$ are irreducible, we must merely check that $V(f_{i})$ does not have any component contained in $X_{1} \setminus (X_{1}\cap X_{2})$.  But $X_{1}\setminus (X_{1}\cap X_{2})$ is just $D_{0}$, and $D_{i} \cap D_{0} = \emptyset$ for $i > 0$.

As the $f_{i}$'s have unbounded degree and are irreducible, their zero sets could not be the fibers of a morphism (or even a rational map!) to a curve.
\qed\\
To prove Theorem \ref{DDArb}, we choose algebraically independent $y, z$ in the ring of regular functions on $X$ for which $V(y)$ and $V(z)$ are irreducible, and construct $f_{i}$ and $D_{i}$ as above (replacing the $D_{i}$ with its finite set of connected components at each stage, if necessary).  Any function to a curve with the $D_{i}$ as fibers would factor through the map to $\mathbb{A}^{2}$ given by $y$ and $z$.
%% Question: does this automatically imply the $D_{i}$ are irreducible?

\rem{This procedure is by no means unique.  For instance, one could replace $f_0$  and $f_1$ by any other two irreducible polynomials with no common roots in $\mathbb A^2_k$, and the above construction works.
However, the following questions remain open: 
\begin{enumerate}
\item Does there exist an example as above where the curves are all smooth?
\item In any example as above, is the geometric genus of the curves necessarily unbounded? (That is, could we find such a counterexample consisting of only rational curves?)
\item In any example as above, does there necessarily exist a divisor $D$ such that $\#(D\cap D_{i})$, the \textit{set-theoretic intersection}, is unbounded?
\end{enumerate}
}
We now prove Theorem \ref{DDAU}, that any \textit{uncountable} set of disjoint divisors must form a family.  In this proof, $X$ will be normal and affine; the theorem follows immediately for all normal varieties from this case.
\prop{ Let $X$ be a  normal quasi-projective variety over an {\it uncountable} algebraically closed field $k$.  Let $\{D_{i}\}_{i\in I}$ be an uncountable collection of  pairwise-disjoint, reduced, connected, codimension one closed subvarieties of $X$.

Then there exist a smooth projective curve $C$ defined over $k$ and a dominant morphism $\varphi: X\rightarrow C$ so that for any $i\in I$ the divisor $D_i$ is contained in a  fiber of $f$.

}
\proof{

Let $X\subset \bar X$ be a projective model of $X$. We may assume that $\bar X$ is normal.
Let $\bar D_i$ be the closure of $D_i$ in $\bar X$.  Note that if there is a subset  $I_0\subset I$ with $\#I_0\geq \rho_w(\bar X)+2$,  such that $\bar D_i$, $i\in I_0$ are disjoint, then we can apply  theorem \ref{DDP} for $\bar X$ to get a map $g:\bar X\to C$ such that all $\bar D_i, i\in I_0$, are contained in the fibers of $g$. Since $D_i, i\in I$ are disjoint, we can take $\varphi$  the restriction of $g$ to $X$. 

 There is an alteration $f: \bar Y\to \bar X$, such that $\bar Y$ is smooth and for $Y= f^{-1}(X)$, we have $Y_{\infty}=\bar Y\setminus Y$ is a simple normal crossings divisor \cite{Il}.  Let $F_i=f^{-1}(D_i)$, then $F_i$ are disjoint and cover $Y$.  Let $\bar F_i\subset \bar Y$ be the closure of $F_i$ in $\bar Y$. Since $\bar Y$ is smooth, each $\bar F_i$ gives a class in the Picard group $Pic\,\bar Y$. Since $NS\,\bar Y$ is countable, we obtain that for $J\subset I$ uncountable, the divisors $\bar F_j, j\in J$ have all the same class  $\alpha\in NS\,\bar Y$.

 \lemma{ There is an infinite  subset $J'\subset J$, a finite set of irreducible divisors $(E_t)_{t\in T}\subset Y_{\infty}$ and $M>0$ such that for any $i,j\in J'$ one has $\bar F_i\cdot \bar F_j=\sum\limits_{t\in T} a_tE_t$ with $0\leq a_t\leq M$. }
 \proof{First note that since the divisors $(\bar F_j)_{j\in J}$ intersect only on $Y_{\infty}$, for a fixed $j\in J$, any intersection $\bar F_{j}\cap \bar F_{j'}$ with $j'\in J$ is supported on components of  $\bar F_{j}\cap Y_{\infty}$, also the intersections $\bar F_{j}\cap \bar F_{j'}$ have all the same class $\alpha^2$.
 
 Consider $j_0\in J$. Since $J$ is uncountable, there is an uncountable subset $J_1\subset J$ such that for any $j, j'\in J_1$ one has $H_0:=\bar F_{j_0}\cdot \bar F_{j}=\bar F_{j_0}\cdot \bar F_{j'}$ as a divisor (not only a class) on $Y_{\infty},$ one may also assume that this intersection is nonzero. Fix now $j_1\in J_1$. Similarly, one finds an uncountable subset $J_2\subset J_1$ such that  for any $j, j'\in J_2$ one has $H_1:=\bar F_{j_1}\cdot \bar F_{j}=\bar F_{j_1}\cdot \bar F_{j'}$.

By the same procedure, we construct inductively the sets $J_n$ and the divisors $H_n$. Since all the divisors $\bar F_j$ have the same class, after a finite number of steps we should obtain $H_{n+r}\subset \cup_{m< n} H_m$ for all $r\geq 0$.  Then $J'=J_n$ and $T$ the set of irreducible components of $\cup_{m< n} H_m$ works. \\ 
 
  \qed\\}
 %Fix $j_0\in J$. Then all the intersections $\bar F_{j_0}\cap \bar F_j$, $j\in J$ are supported on $\bar F_{j_0}\cap Y_{\infty}$, hence  $\bar F_{j_0}\cap \bar  F_j$ is a sum of components of $\bar F_{j_0}\cap Y_{\infty}$ with some multiplicities.  Since $J$ is uncountable, one can find $J'\subset J$ uncountable, such that for any $j\in J'$, the intersection   $\bar F_{j_0}\cap \bar F_j$  on $Y_{\infty}$ is the same divisor on $Y_{\infty}$.  
 
 The lemma above gives the following bound on $\bar X$: there is an integer $N$ such that for any closed $Z\subset \bar X$ that is a (set-theoretic) component of the intersection of $\bar D_{j}$ and $\bar D_{j'}, j\in J'$ we have that locally in $\mathcal O_{X, Z}$,   the ideal of the intersection of $\bar D_{j}$ and $\bar D_{j'}$ is contained in at most the $N^{th}$ power of the maximal ideal  $\mathfrak m_{Z}^N\subset \mathcal O_{X, Z}$ of $Z$. Hence after a finite number of blow-ups $\tilde X\to \bar X$ centered at $X_{\infty}$,  the strict transforms $\tilde D_j$ of $\bar D_j$  do not intersect. Now we can apply theorem \ref{DDP} to $\tilde X$ and the family $(\tilde D_j)_{j\in J'}$ to get a map $f:\tilde X\to C$, such  that the restriction $\phi$ of $f$ to $X$ satisfies the conclusion of the theorem.

 \qed\\}

\paragraph{\bf Acknowledgements.}
The first author acknowledges that the article was prepared within the framework of a subsidy granted to the HSE by the Government of the Russian Federation for the implementation of the Global Competitiveness Program.  The first author was also supported by a Simons Travel Grant. 
The second author thanks the University of Chicago and the Embassy of France in United States  for their support for a short-term visit  to the University of Chicago in April, 2015. The third author thanks Madhav Nori, Keerthi Madapusi-Pera, Mihnea Popa, and Jesse Kass for helpful conversations, and gratefully acknowledges the support of NSF Grant DMS-1400683.  We gratefully acknowledge Burt Totaro and John Christian Ottem for pointing us in the direction of \cite{TotaroDisjDiv} and \cite{PereiraFol}, after this paper was first posted on the ArXiv.  This paper was conceived at Workshop \#1444 at Oberwolfach, and we gratefully acknowledge the hospitality of the Matematisches Forschungsinstitut Oberwolfach.

\vspace{1 cm}

\end{document}